\documentclass[english]{article}
\usepackage[latin9]{inputenc}
\usepackage{amstext}
\usepackage{wasysym}

\makeatletter

\date{}

\makeatother

\usepackage{babel}
\begin{document}

\title{A Note on Intervals in the Hales-Jewett Theorem}

\author{Imre Leader and Eero R\"aty \footnote{Centre for Mathematical Sciences,
Wilberforce Road,
Cambridge CB3 0WB,
UK,  {\tt i.leader@dpmms.cam.ac.uk} and 
{\tt epjr2@cam.ac.uk}}}

\maketitle
\begin{abstract}
The Hales-Jewett theorem for alphabet of size 3 states that whenever
the Hales-Jewett cube $[3]^{n}$ is $r$-coloured there is a monochromatic
line (for n large). Conlon and Kamcev conjectured that, for any $n$,
there is a 2-colouring of $[3]^{n}$ for which there is no monochromatic
line whose active coordinate set is an interval. In this note we disprove
this conjecture.
\end{abstract}

\section{Introduction}

In order to state the Hales-Jewett theorem we need some notation.
Given positive integers $k$ and $n$ let $[k]^{n}$ be the set of
all words in symbols $\{1,\dots,k\}$ of length $n$. A set $L \subset 
[k]^{n}$ is called
a combinatorial line if there exist a nonempty set $S\subset [n]$
and integers $a_{i}\in[k]$ for all $i\not\in S$ such that $L=\{(x_{1},\dots,x_{n})\,:\,x_{i}=a_{i}\, \mbox{for all}\, i\not\in S,\,x_{i}=x_{j}\, 
\mbox{for all}\, i,j\in S\}$.
The set $S$ is called the active coordinate set of $L$. \\

\textbf{Theorem (Hales-Jewett \cite{key-2}). }For any $k$ and $r$ there exists
$n$ such that whenever $[k]^{n}$ is $r$-coloured 
there is a monochromatic combinatorial line. \\

As noted by Conlon and Kamcev in \cite{key-1}, by following Shelah's
proof \cite{key-3} of the Hales-Jewett theorem it can be shown
that one can always find a monochromatic combinatorial line whose
active coordinate set $S$ is an union of at most $HJ(k-1,r)$ intervals,
where $HJ(k-1,r)$ is the smallest integer $n$ for which the Hales-Jewett
theorem holds for $k-1$ and $r$. 

In the case $k=3$, since $HJ(2,r)=r$, this
says that one can always find a monochromatic line whose active coordinate
set is a union of at most $r$ intervals.
Conlon and Kamcev proved in \cite{key-1}
that this bound is tight for $r$ odd: in other words, they showed that for
each odd $r$ there is an $r$-colouring of $[3]^{n}$ (for any $n$) for which every 
monochromatic line has active coordinate set made up of at least $r$ intervals.
They conjectured
that this would also be the case for $r$ even. In particular, for $r=2$, they
conjectured that for all $n$ there exists a 2-colouring
of $[3]^{n}$ for which there exists no monochromatic combinatorial
line whose active coordinate set is an interval. 

In this note we will
prove that, perhaps surprisingly, their conjecture is false. This can be stated
in the following form. \\

\textbf{Theorem 1.} For all sufficiently large $n$, whenever $[3]^{n}$
is 2-coloured there exists a monochromatic combinatorial line 
whose active coordinate set is an interval.

\section{The proof of Theorem 1}

The idea of the proof is as follows. By applying Ramsey's theorem,
we will pass to a subspace on which the colour of a word depends only
on its `pattern' of intervals. Once this is done, we can consider
some particular small patterns.

For a word $w$ let $\overline{w}$ be obtained from $w$ by contracting
every interval on which $w$ is a constant to a single letter. We will
consider the particular words $s_{1}=132$ , $s_{2}=1232$, $s_{3}=1312$,
$s_{4}=13232$ and $s_{5}=13132$. Set $t_{i}$ to be the length of
the word $s_{i}$. Put $n_{0}=4$ and for $1 \leq i \leq 5$
let $n_{i}=R^{(t_{i}-1)}(n_{i-1})$ where $R^{(t)}(s)=R^{(t)}(s,s)$
is the $t$-set Ramsey-number. Finally set $n=n_{5}+1$. 

Let $c$ be any 2-colouring of $[3]^{n}$.
For a word $w$ we define the set of \textit{breakpoints} $T(w)$
by $T(w)=\{a_{1},\dots,a_{m}\}$ if $w_{a_{i-1}\text{+1}}=\dots=w_{a_{i}}$
and $w_{a_{i}}\neq w_{a_{i}+1}$ for all $1\leq i\leq m+1$, with
the convention $a_{0}=0$, $a_{m+1}=n$. For example, $w=1122333111$
has breakpoints $T(w)=\{2,4,7\}$.

Let $s$ be a sequence of length $t$ and $T_{1}=\{a_{1},\dots,a_{m}\}\subset [n-1]$
with $|T_{1}|=t-1$. We say that $w\in[3]^{n}$ has breakpoints in
$T_{1}$ with \textit{pattern} $s$ if $T(w)=T_{1}$ and $\overline{w}=s$.
For example $w=1122333111$ has breakpoints $T(w)=\{2,4,7\}$ with
pattern $s=1231$. Note that if $\overline{w}=s$ then there exists
a unique set $T_{1}$ of size $|s|-1$ for which $w$ has breakpoints
$T_{1}$ with pattern $s$.

Set $T_{5}=[n-1]$. Suppose that $|T_{i}|\geq n_{i}$ is given, and
recall that $t_{i}$ is the length of the word $s_{i}$ defined at
the start of the proof. For all $A\in[n-1]^{(t_{i}-1)}$ define
$w^{A}$ to be the unique sequence which has breakpoints $A$ with
pattern $s_{i}$. 

Now $c$ induces a 2-colouring $c_{i}$ on the set $T_{i}^{(t_{i}-1)}$
given by $c_{i}(A)=c(w^{A})$. Hence by Ramsey's theorem and the choice
of $n_{i}$'s it follows that there exists $T_{i-1}\subset T_{i}$
with $|T_{i-1}|\geq n_{i-1}$ such that $T_{i-1}^{(t_{i}-1)}$ is
monochromatic for the colouring $c_{i}$ , say with colour $d_{i}$.

Thus we obtain sets $T_{0}\subset T_{1}\subset\dots\subset T_{5}$
with $|T_{0}|\geq4$ and colours $d_{1}\dots,d_{5}$ such that $c_{i}$
restricted to $T_{i}^{(t_{i+1}-1)}$ is constant with value $d_{i+1}$. Note
that it is impossible to choose colours $d_{1},\dots,d_{5}$ without
at least one of the following sets 
\[
\begin{array}{c}
N_{1}=\{d_{1},d_{2}\}\\
N_{2}=\{d_{1},d_{3}\}\\
N_{3}=\{d_{2},d_{4}\}\\
N_{4}=\{d_{3},d_{5}\}\\
N_{5}=\{d_{1},d_{4},d_{5}\}
\end{array}
\]
having just one element (i.e. all colours being equal). Indeed, if
all the sets $N_{1},N_{2},N_{3},N_{4}$ contain both colours we must
have $d_{2}=d_{3}$ and $d_{1}=d_{4}=d_{5}$, which implies that $|N_{5}|=1$. 

Let $a_{1}<a_{2}<a_{3}<a_{4}$ be elements of $T_{0}$. We will use
the shorthand $w=[b_{1}b_{2}b_{3}b_{4}b_{5}]$ for the word which
has $w_{i}=b_{j}$ for all $a_{j-1}<i\leq a_{j}$ with the convention
$a_{0}=0$ and $a_{5}=n$. Note that we will allow $b_{i}=b_{i+1}$.
Hence $T(w)\subset\{a_{1},\dots,a_{4}\}\subset T_{0}$ and $\overline{w}=\overline{b_{1}b_{2}b_{3}b_{4}b_{5}}$.

Set $w_{1}=[13332]$, $w_{2}=[12232]$, $w_{3}=[13112]$, $w_{4}=[13232]$
and $w_{5}=[13132]$. It is easy to verify that for all $i$ we have
$\overline{w_{i}}=s_{i}$ and also $T(w_{i})\subset T_{0}\subset T_{i-1}$.
Furthermore set $v_{1}=[11132]$, $v_{2}=[11232]$, $v_{3}=[13122]$
and $u_{1}=[13222]$. As before it is easy to verify that $\overline{v_{i}}=s_{i}$,
$\overline{u_{1}}=s_{1}$, and by construction $T(v_{i})\subset T_{0}\subset T_{i-1}$
and $T(u_{1})\subset T_{0}$. Thus by definition of the sets $T_{i-1}$
it follows that $c(w_{i})=d_{i}$, $c(v_{i})=d_{i}$ and $c(u_{1})=d_{1}$. 

It is straightforward to verify that 
\begin{itemize}
\item $v_{1},w_{2},w_{1}$ forms a combinatorial line $L_{1}$ with $S_{1}=\{a_{1}+1,\dots,a_{3}\}$
\item $w_{3},u_{1},w_{1}$ forms a combinatorial line $L_{2}$ with $S_{2}=\{a_{2}+1,\dots,a_{4}\}$
\item $v_{2},w_{2},w_{4}$ forms a combinatorial line $L_{3}$ with $S_{3}=\{a_{1}+1,\dots,a_{2}\}$
\item $w_{3},v_{3},w_{5}$ forms a combinatorial line $L_{4}$ with $S_{4}=\{a_{3}+1,\dots,a_{4}\}$
\item $w_{5},w_{4},w_{1}$ forms a combinatorial line $L_{5}$ with $S_{5}=\{a_{2}+1,\dots,a_{3}\}$
\end{itemize}
It is clear that the colours used to colour elements of the line $L_{i}$
are exactly the colours in the set $N_{i}$. As observed earlier, one
of the sets $N_{i}$ contains only one colour, which implies that
the associated line $L_{i}$ is monochromatic. Since all the sets
$S_{i}$ are intervals, this completes the proof.

\end{document}